\documentclass[12pt,twoside,a4paper]{article}
\usepackage{stmaryrd}
\usepackage{amsfonts}
\usepackage{amsmath,amssymb}
\usepackage{mathrsfs}
\usepackage{hyperref}

\newtheorem{thm}{Theorem}[section]

\newtheorem{lem}[thm]{Lemma}
\newtheorem{prop}[thm]{Proposition}

\numberwithin{equation}{section}\allowdisplaybreaks

\def\le{\leqslant}
\def\ge{\geqslant}
\def\leq{\leqslant}
\def\geq{\geqslant}

\begin{document}

\begin{center}
\Large\bf  Ill-posedness for the Navier-Stokes equations in critical Besov spaces $\dot B^{-1}_{\infty,q}$

 \vspace*{0.5cm}\normalsize
   Baoxiang Wang\footnote{Email:
  wbx@math.pku.edu.cn}  \  

\vspace*{0.5cm} {\it \footnotesize LMAM, School of Mathematical Sciences, Peking University,
Beijing 100871   China  }\\

\vspace*{0.5cm}

\begin{minipage}{13.5cm}
\footnotesize \bf Abstract. \rm  We study the Cauchy problem for the incompressible
Navier-Stokes equations in two and higher spatial dimensions
\begin{align}\label{NS}
u_t -\Delta u+u\cdot \nabla u +\nabla p=0, \ \  {\rm div} u=0, \ \  u(0,x)= \delta u_0.
\end{align}
For arbitrarily small $\delta>0$, we show that the solution map $\delta u_0 \to u$ in critical Besov spaces $\dot B^{-1}_{\infty,q}$ ($\forall \ q\in [1,2]$) is discontinuous at origin. It is known that the Navier-Stokes equation is globally well-posed for small data in $BMO^{-1}$ (\cite{KoTa01}). Taking notice of the embedding $\dot B^{-1}_{\infty,q} \subset BMO^{-1}$ ($q\le 2$), we see that for sufficiently small $\delta>0$, $u_0\in  \dot B^{-1}_{\infty,q} $ ($q\le 2$) can guarantee that \eqref{NS} has a unique global solution in $BMO^{-1}$, however, this solution is instable in $ \dot B^{-1}_{\infty,q} $ and the solution can have an inflation in $\dot B^{-1}_{\infty,q} $ for certain initial data.  So, our result indicates that two different topological structures in the same space may determine  the well and ill posedness, respectively. \\

\\

\end{minipage}
\end{center}

\section{Introduction} \label{sect1}

\medskip

We study the ill-posedness for the Cauchy problem of the  incompressible
Navier-Stokes equations (NS):
\begin{align}\label{NS1}
u_t -\Delta u+u\cdot \nabla u +\nabla p=0, \ \  {\rm div} u=0, \ \  u(0,x)= \delta u_0.
\end{align}
where $t\in \mathbb{R}^{+}=[0,\infty)$, $x\in \mathbb{R}^{n} \ (n\geq 2)$,   $u=(u_1,...,u_n)$ denotes the flow
velocity vector and $p(t,x)$ describes the scalar pressure. $0<\delta \ll 1$ denotes a small parameter, $\nabla=(\partial_1,...,\partial_n)$ is the gradient operator, $\Delta=\partial^2_1+...+\partial^2_n$ is the Laplacian,
$u_0(x)=(u^0_1,...,u^0_n)$ is a given velocity with $\mathrm{div}
u_0=\partial_1 u^0_1+...+\partial_n u^0_n =0$.
It is easy to see that \eqref{NS1} can be rewritten as the following
equivalent form:
\begin{align}
  u_t -\Delta u+\mathbb{P}\ \textrm{div}(u\otimes u)=0, \ \
u(0, x)= \delta u_0,
 \label{2.1}
\end{align}
where $\mathbb{P}=I-\nabla \Delta^{-1}{\rm div}$ is the
projection operator onto the divergence free vector fields.

It is known that \eqref{NS1} is essentially equivalent to the following integral equation:
\begin{align} \label{1.3}
u(t) =   e^{ t \Delta } u_0 +  \int^t_0  e^{ (t-\tau) \Delta } \mathbb{P}\ \textrm{div}(u\otimes u)(\tau) d\tau.
\end{align}

Note that \eqref{NS1} is scaling invariant in the following sense:
if $u$ solves \eqref{NS1},  so does   $u_{\lambda}(t, x)=\lambda u(\lambda^{2 }t, \lambda x)$ and $p_\lambda (t, x)=\lambda^2 p(\lambda^{2 }t, \lambda x)$ with initial data $\lambda u_0(\lambda x)$.
 A function space $X$ defined in $\mathbb{R}^n$ is said to be a \textit{ critical space} for \eqref{NS1}
if the  norms of $u_{\lambda}(0,x)$ in $X$ are equivalent for
all $\lambda>0$ (i.e., $\|u_\lambda(0,\cdot)\|_X\sim \|u_0\|_X$). It is easy to see that the following spaces are critical spaces for NS:
\begin{align}
& \mbox{Besov  type spaces:} \     \dot{B}^{n/p-1}_{p, q} \subset     \dot{B}^{-1}_{\infty, r}  \subset     \dot{B}^{-1}_{\infty, \infty}\ (p< \infty, \ q\le r), \nonumber\\
& \mbox{Triebel  type spaces:} \  L^n = \dot{F}^{0}_{n, 2} \subset   \dot{F}^{n/p-1}_{p, q}    \subset \dot{F}^{-1}_{\infty, r} \ (n< p< \infty). \nonumber
\end{align}

NS \eqref{NS1} has been extensively studied  in the past twenty years; cf. \cite{BaBiTa12, Can97,  Ch99, ChGaPa11, EsSeSv03, FoTe89, Iw10, Ka84, KeKo11, KoTa01, Pl96} and references therein. For the initial data in critical Besov type spaces, Cannone \cite{Can97}, Planchon \cite{Pl96}  and  Chemin \cite{Ch99} obtained global solutions in 3D for small data  in critical Besov spaces $\dot{B}^{3/p-1}_{p, q}$ for all $p<\infty, \ q\le \infty$. In the case $p=\infty$, Bae, Biswas and  Tadmor \cite{BaBiTa12} can show the  global well posedness of the solutions in 3D for the sufficiently small initial data in $\dot{B}^{-1}_{\infty, q} \cap \dot{B}^{0}_{3, \infty}$ with $1\le q<\infty$.
Bourgain and Pavlovic \cite{BoPa08} showed the  \textit{ill-posedness} of NS in $\dot{B}^{-1}_{\infty,\infty}$, i.e., the solution map is discontinuous in $\dot{B}^{-1}_{\infty,\infty}$,   Germain \cite{Ge08} proved that the solution map of NS is not $C^2$ in $\dot{B}^{-1}_{\infty, q} $ for any $q>2$,  Yoneda \cite{Yo10} can show that the solution map is discontinuous in $\dot{B}^{-1}_{\infty, q} $ for any $q>2$ and in fact, he constructed a logarithmic type Besov space $V$ very near to  $\dot{B}^{-1}_{\infty, 2}$ so that NS is ill-posed in $V$.  Up to now, the largest Besov-type space on initial data for which NS is well-posed is still unknown.

For the initial data in critical Triebel type spaces,  Koch and Tataru  \cite{KoTa01} obtained the global well-posedness for small initial data in $BMO^{-1}= \dot{F}^{-1}_{\infty, 2}$, Yoneda \cite{Yo10} pointed out that his argument implies that NS is ill-posed in $\dot{F}^{-1}_{\infty, q}$ for all $q>2$ (cf. also \cite{DeYa13}). Now let us recall Koch and Tataru's result:

\begin{thm}\label{NSWellposed} {\rm (\cite{KoTa01,Iw13}, Global well-posedness for small data in $BMO^{-1}$)}
Let   $u_0 \in {BMO}^{-1}$ with $\|u_0\|_{BMO^{-1}} \le 1$. Then there exists a $\delta_0>0$ such that for any $0<\delta \le \delta_0$,  NS has a unique global solution    $u(\delta,t) \in X$  with $\|u\|_{X} \le C\delta$, where
$$
\|u\|_X := \sup_{t>0} t^{1/2}\|u(t,\cdot)\|_{L^\infty(\mathbb{R}^n)} + \sup_{x\in \mathbb{R}^n, R>0} |B(x, R)|^{-1/2} \|u(t,y)\|_{L^2_{t,y}([0,R^2]\times B(x,R))}.
$$
\end{thm}
Moreover, Auscher, Dubois and Tchamitchian \cite{AuDuTc06}, Iwabuchi and Nakamura \cite{Iw13} obtained that Koch and Tataru's solution is stable and belongs to $L^\infty (0,\infty; BMO^{-1})$.
Recalling that the inclusions $ \dot{B}^{s}_{p, \min\{p,q\}} \subset \dot{F}^{s}_{p, q} \subset \dot{B}^{s}_{p, \max\{p,q\}}$ ($q>1$),  one sees that  $\dot{B}^{-1}_{\infty, 2}\subset \dot{F}^{-1}_{\infty, 2}$.
According to their results, $u_0\in \dot{B}^{-1}_{\infty, 2}$ and $0<\delta\ll 1$ imply that NS has a unique global solution in $BMO^{-1}$. So,  it seems natural to conjecture that NS is also globally well-posed in $\dot{B}^{-1}_{\infty, 2}$ for sufficiently  small $\delta>0$.
However, in this paper we will show the following negative result:

\begin{thm}\label{NSIll} {\rm (Ill-posedness of the solution)}
Let $n\ge 2$,   $1\le q\le 2$, $0<\delta \ll 1$. Let $u(\delta, t)$ be Koch and Tataru's solution of \eqref{1.3}.  Then the solution map $\delta u_0 \to u(\delta, t)$ in $\dot{B}^{-1}_{\infty, q}(\mathbb{R}^n)$ is discontinuous  at $\delta=0$. More precisely, for any $N\gg 1$,  there exists $u_0 \in \dot{B}^{-1}_{\infty, q}(\mathbb{R}^n)$ with $\|u_0\|_{\dot{B}^{-1}_{\infty, q}} \lesssim 1$, $t  \le 1/N$ such that
$$\left\|   u (\delta, t )\right\|_{\dot{B}^{-1}_{\infty, q}(\mathbb{R}^n)}  \ge (\log N)^{1/2 q}.$$
\end{thm}
We remark that our result also holds for all $2\le q<\infty$. In the proof of Theorem \ref{NSIll} we have no any conditions on $q\ge 1$. The techniques used in this paper are quite different from the arguments as in \cite{BoPa08, Ge08, Yo10}.

Throughout this paper, $C\ge 1, \ c\le 1$ will denote constants which can be different at different places, we will use $A\lesssim B$ to denote   $A\leqslant CB$. We denote    by $L^p=L^p(\mathbb{R}^n)$ the Lebesgue space on which the norm is written as $\|\cdot\|_p$.
Now let us recall the definition of Besov type spaces; cf. \cite{BL, Tr}. Let $\psi: \mathbb{R}^n\rightarrow
[0, 1]$ be a smooth cut-off function which equals $1$ on the ball $B(0,5/4):=\{\xi\in \mathbb{R}^n: |\xi|\le 5/4\}$ and equals $0$ outside the ball $B(0,3/2)$.
Write
\begin{align} \label{phi}
\varphi(\xi):=\psi(\xi)-\psi(2\xi),  \ \ \varphi_j(\xi)=\varphi(2^{-j}\xi),
\end{align}
$\triangle_j:=\mathscr{F}^{-1}\varphi_j \mathscr{F}, \  j\in \mathbb{Z}$
are said to be the dyadic decomposition operators\footnote{$\triangle$ and $\Delta$ are different notations in this paper.}. One easily sees that ${\rm supp} \varphi_j \subset B(0, 2^{j+1})\setminus B(0, 2^{j-1})$ and
\begin{align} \label{phij}
  \varphi_j(\xi)=1, \ \ \mbox{if }  \xi \in B(0,5 \cdot 2^{j-2})\setminus B(0,3\cdot 2^{j-2}).
\end{align}
Let $s\in \mathbb{R}$, $1\le p,q\le \infty$.  The norms in homogeneous Besov and Triebel spaces are defined as follows:
\begin{align}
\label{Besov} \|f\|_{\dot{B}^s_{p, q}}=
 \left(\sum_{j=-\infty}^{+\infty}2^{jsq}\|\triangle_j f\|^q_{p}\right)^{1/q},  \ \
\|f\|_{ \dot{F}^s_{p, q} } = \left\|
\left(\sum_{j=-\infty}^{+\infty} |2^{js } \triangle_j
f |^q \right)^{1/q} \right\|_{p}
\end{align}
with the usual modification for $q=\infty$, where we further assume $p\neq \infty$ in $\dot{F}^s_{p, q}$. Since the definition of $\dot{F}^s_{\infty, q}$ is slightly different from \eqref{Besov},   we leave it into Section \ref{other}. In this paper, we will frequently us the following Bernstein's multiplier estimate (see \cite{BL, WaHuHaGu11}):
\begin{prop}\label{Bern} {\rm (Bernstein's multiplier estimate)} Let $L\in \mathbb{N}, \  L>n/2, \ \theta=n/2L$. We have
\begin{align}
\label{1.7}
\|\mathscr{F}^{-1} \rho\|_{1} \lesssim \|\rho\|^{1-\theta}_{2} \left(\sum^n_{i=1} \|\partial^L_i \rho\|_2^{\theta} \right).
\end{align}
\end{prop}
Let us consider Taylor's expansion of Koch and Tataru's solution $u(\delta,t )$ of NS, one can find a small $\delta_0>0$ such that for any $0<\delta \le \delta_0$,
\begin{align}
\label{1.8}
u(\delta, t) = \sum^\infty_{r=0} \frac{\delta^r}{r!} \frac{\partial^r u}{\partial \delta^r}(0,t) \ \ \ \mbox{in} \ \ X \cap L^\infty(0,\infty; BMO^{-1}).
\end{align}
Recall that Koch and Tataru applied the contraction mapping method to show the global well-posedness of NS for small data in $BMO^{-1}$, we see that \eqref{1.8} can be obtained by iterations.
If NS is globally well-posed in $\dot{B}^{-1}_{\infty, q}$ ($q\in [1,2]$), the solution must coincide with the solution $u(\delta,t)$ in $ X$. So, if we can find some initial data $u_0 \in  \dot{B}^{-1}_{\infty, q}$ such that the second iteration $ \delta^2 \partial_\delta^2  u (0,t)/2$ has an inflation in  $\dot{B}^{-1}_{\infty, q}$ and the norm of the other terms in \eqref{1.8} are much less than that of the second iteration, then our result is shown. More precisely, for some subset $\mathbb{A} \subset \mathbb{N}$,
\begin{align}
\| u(\delta, t)\|_{\dot{B}^{-1}_{\infty, q}}  \ge & \left( \sum_{j\in \mathbb{A}} 2^{-jq} \|\triangle_j  u(\delta, t) \|^q_\infty   \right)^{1/q} \nonumber\\
 \ge & \frac{\delta^2}{2} \left( \sum_{j\in \mathbb{A}} 2^{-jq} \|\triangle_j \partial_\delta^2  u (0,t) \|^q_\infty   \right)^{1/q} \label{1.10} \\
 & - \left( \sum_{j\in \mathbb{A}} 2^{-jq} \left\|\triangle_j \sum_{r\ge 1, \ r\neq 2} \frac{\delta^r}{r!} \partial_\delta^r  u (0,t) \right\|^q_\infty   \right)^{1/q}. \label{1.11}
\end{align}
We will show that  \eqref{1.10}  contributes the main part and it will be quite large after any short time and  \eqref{1.11} is much less than  \eqref{1.10}.

The paper is organized as follows. In Section \ref{sect2} we prove Theorem \ref{NSIll} for higher dimensions $n\ge 3$ and in Section \ref{sect3} we continue its proof for 2D case. In Section \ref{other} we compare the solutions in different critical spaces.

\section{Proof of Theorem \ref{NSIll}: $n\ge 3$} \label{sect2}

\subsection{Estimates on second iteration }

Let $u(\delta,t)$ be the solution of \eqref{NS1} with $0<\delta \ll 1$,  we see that
\begin{align}
& \label{seconditerat} u(\delta, t)|_{\delta=0} =0, \ \    \frac{\partial  u}{\partial\delta}\Big|_{\delta=0} = e^{t\Delta } u_0, \nonumber\\
& \frac{\partial^2  u}{\partial\delta^2}\Big|_{\delta=0} = 2 \int^t_0  e^{(t-\tau)\Delta }  \mathbb{P} (e^{\tau\Delta } u_0 \cdot \nabla ) e^{\tau\Delta } u_0 d\tau.
\end{align}
Our aim is to choose a suitable $u_0 \in \dot B^{-1}_{\infty, q}$ which is localized in one dyadic frequency $|\xi| \sim 2^k$ and to find a time $t\sim 2^{-2k}$ verifying
 $ \| \frac{\partial^2 u}{\partial\delta^2}(0,t) \|_{\dot{B}^{-1}_{\infty, q}(\mathbb{R}^n)}  \gtrsim k^{1/q}.$
Let $k\in 16\mathbb{N}=\{16, 32, 48,...\}$, $k\gg 1$, $0<\varepsilon \ll 1$\footnote{$\varepsilon>0$ will be chosen below as in \eqref{2.78}, $k$ can be arbitrarily large.}. Denote
\begin{align}
a_l  = 2^l (\varepsilon, 2\varepsilon, \sqrt{1-5\varepsilon^2},0,...,0), \ \
b_l =a_l/2, \ \ c_k   = \frac{2^k}{\sqrt{n}} (1,1,...,1).  \label{ck}
\end{align}
We write
\begin{align}
\mathbb{N}_k  =  \{l \in 4\mathbb{N}: \ k/4\le l \le k/2\}.  \label{ck}
\end{align}
Let $\psi$ be as in \eqref{phi} and $\varrho(\xi) = \psi(4 \xi)$ and
\begin{equation}\label{Phi}
        \begin{cases}
                 \widehat{\Phi^+_l}(\xi) = e^{{\rm i}\xi a_l} (\varrho (\xi-c_k-b_l) + \varrho (\xi-c_k+b_l)),\\
                 \widehat{\Phi^-_l}(\xi) = e^{{\rm i}\xi a_l} (\varrho (\xi+c_k-b_l) + \varrho (\xi+c_k+b_l)).\\
        \end{cases}
\end{equation}
We now introduce the initial data $u_0=(u^0_1,...,u^0_n)$:
\begin{equation}\label{u0}
        \begin{cases}
                 \widehat{u^0_1} (\xi) & = 2^k \sum_{l\in \mathbb{N}_k}  ( \widehat{\Phi^+_l}(\xi) + \widehat{\Phi^-_l}(\xi)),\\
                 \widehat{u^0_2} (\xi) &  = -\frac{\xi_1}{\xi_2} \widehat{u^0_1} (\xi) = - 2^k \sum_{l\in \mathbb{N}_k} \frac{\xi_1}{\xi_2} ( \widehat{\Phi^+_l}(\xi) + \widehat{\Phi^-_l}(\xi))
    \end{cases}
\end{equation}
and $u^0_3(x) = ... = u^0_n (x) =0.$  One can rewrite  $ u^0_1 $ and $ u^0_2 $  as
\begin{align}
 u^0_1 (x) & = 2^k \sum_{l\in \mathbb{N}_k} ( \cos{ (x + a_l)(c_k+b_l)} + \cos{ (x + a_l)( b_l - c_k)})    \check{\varrho}  (x+a_l)
   \label{u01r}\\
 u^0_2  (x) &      = -u^0_1(x) + 2^k \sum_{l\in \mathbb{N}_k} \mathscr{F}^{-1} \frac{\xi_2-\xi_1}{\xi_2} ( \widehat{\Phi^+_l}(\xi) + \widehat{\Phi^-_l}(\xi)).
   \label{u02r}
  \end{align}
Now we give some explanations to the initial data. We can assume that $\varrho(\xi) $ is radial so that $\check{\varrho}$ is a real function.  From $\widehat{u^0_2}= - (\xi_1/\xi_2) \widehat{u^0_1}$ it follows that  ${\rm div} u_0=0$. Introducing  $c_k$  is to guarantee that $\widehat{u^0_1}$ and $\widehat{u^0_2}$ are supported in the dyadic $\{\xi: 2^{k-1} \le |\xi| \le 2^{k+1}\}$, it follows that the frequency of $u_0$ is very high if $k$ is very large. The nonlinear interaction will pullback $\widehat{u^0_1}*\widehat{u^0_1}$ into much lower frequencies. In order to $\widehat{u^0_1}*\widehat{u^0_1}$ concentrated in different dyadic regions as many as possible, different $b_l$  is introduced. $a_l$ is used for controlling the superpositions of all $\check{\varrho}(\cdot+a_l)$ in \eqref{u01r}.

\begin{lem}\label{NSlem1}
Let   $1\leq q\le  \infty$, $k\gg -\log \varepsilon$. Then
\begin{align} \label{boundu012}
\|u^0_i\|_{  \dot{B}^{-1}_{\infty,
q} }\lesssim 1, \ \ i=1,2.
\end{align}
\end{lem}
{\it Proof.} Since $\check{\varrho}$ is a Schwartz function, we have
\begin{align}
|\check{\varrho} (x)| \lesssim (1+|x|)^{-N}, \ N \gg 1.
   \label{rapiddecay}
  \end{align}
It follows from \eqref{u01r} and \eqref{rapiddecay} that
\begin{align}
\| u^0_1  \|_\infty  \lesssim 2^k \left\| \sum_{l\in \mathbb{N}_k} ( 1+|x+a_l|)^{-N} \right\|_\infty  \lesssim 2^k.
   \label{2.11}
  \end{align}
Hence, noticing that ${\rm supp} \widehat{u^0_1}$ is  included  in $\{\xi: |\xi| \sim 2^k\}$, we have
\begin{align}
\| u^0_1  \|_{\dot B^{-1}_{\infty, q}}  \lesssim 1.
   \label{boundu01}
  \end{align}
Applying Young's inequality and Bernstein's multiplier estimate, we have for any $l\in \mathbb{N}_k$ ($k\gg -\log \varepsilon$),
\begin{align}
\left\| \mathscr{F}^{-1} \frac{\xi_2-\xi_1}{\xi_2} e^{{\rm i}\xi a_l} \varrho(\xi+c_k - b_l)  \right\|_\infty & \le   \left\| \mathscr{F}^{-1} \frac{\xi_2-\xi_1}{\xi_2}   \psi(\xi+c_k - b_l)\right\|_1 \|\check{\varrho}\|_\infty \nonumber\\
 & \lesssim  \left\| \mathscr{F}^{-1} \frac{\xi_2-\xi_1+ \varepsilon 2^{l-1}}{\xi_2- 2^k/\sqrt{n} +2\varepsilon 2^{l-1}}   \psi(\xi)\right\|_1 \nonumber\\
 & \lesssim \varepsilon 2^{-k/2}. \label{2.13}
  \end{align}
Similarly, one can estimate the other terms in $\Phi^{\pm}_l$ and we have
\begin{align}
\left\| \mathscr{F}^{-1} \frac{\xi_2-\xi_1}{\xi_2}   \widehat{\Phi^{\lambda}_l}(\xi)\right\|_\infty    \lesssim \varepsilon 2^{-k/2}, \quad \lambda=+,-. \label{2.14}
  \end{align}
Collecting \eqref{u02r}, \eqref{boundu01}--\eqref{2.14},
\begin{align}
\| u^0_2  \|_{\dot B^{-1}_{\infty, q}}  \lesssim 1 + \sum_{l\in \mathbb{N}_k} \left\| \mathscr{F}^{-1} \frac{\xi_2-\xi_1}{\xi_2}   (\widehat{\Phi^+_l}(\xi)+ \widehat{\Phi^-_l}(\xi)) \right\|_\infty \lesssim 1+  \varepsilon k 2^{-k/2}  \lesssim 1.
  \end{align}
So, we have the desired bounds of $u^0_i, \ i=1,2.$ $\hfill \Box$

Recall that for the solution $u=(u_1(\delta,t),...,u_n(\delta,t))$ with initial data $\delta u_0$,
\begin{align}
\frac{\partial^2  u_1}{\partial\delta^2}\Big|_{\delta=0} & = \int^t_0  e^{(t-\tau)\Delta }  \left[ \sum_{i=1}^2 \partial_i (e^{\tau\Delta } u^0_i e^{\tau\Delta } u^0_1) - \partial_1 \sum_{i,j =1}^2 \frac{\partial_i \partial_j}{\Delta} (e^{\tau\Delta } u^0_i e^{\tau\Delta } u^0_j) \right]d\tau \nonumber\\
& := F_1(t,x) -F_2(t,x). \label{secondder}
\end{align}
It follows that
\begin{align}
  & \left(\sum_{j\in \mathbb{N}_k}  2^{-jq} \left\|\triangle_j \frac{\partial^2  u_1}{\partial\delta^2}\Big|_{\delta=0} \right\|^q_\infty \right)^{1/q} \nonumber\\
 & \quad \ge  \left(\sum_{j\in \mathbb{N}_k}  2^{-jq} \left\|\triangle_j  F_1  \right\|^q_\infty \right)^{1/q} -  \left(\sum_{j\in \mathbb{N}_k}  2^{-jq} \left\|\triangle_j  F_2  \right\|^q_\infty \right)^{1/q} \nonumber\\
 & \quad : = A_1-A_2.  \label{secondder}
\end{align}
In the following, our aim is to show that for $t=\varepsilon^2 2^{-2k}$, $A_1 \gtrsim \varepsilon^3 k^{1/q}$ and $A_2  \ll  \varepsilon^3 k^{1/q}$. Since $k$ can be arbitrarily large, one immediately has   $\left\| \partial^2  u_1/\partial\delta^2  |_{\delta=0} \right\|_{\dot B^{-1}_{\infty, q}} \to \infty$ as $k\to \infty.$ For convenience, we denote
\begin{align}
R_i(f,g) & = \int^t_0  e^{(t-\tau)\Delta }  \partial_i  (e^{\tau\Delta } f \  e^{\tau\Delta } g) d\tau; \label{2.18}\\
\widetilde{R}_i(f,g) & = \int^t_0  e^{(t-\tau)\Delta }  \partial_i  (e^{\tau\Delta } f \  e^{\tau\Delta } g - fg) d\tau.  \label{2.19}
\end{align}
It follows that
\begin{align}
R_i(f,g) = \Delta^{-1} ( e^{t\Delta }-1) \partial_i (fg) +  \widetilde{R}_i(f,g).
\end{align}
Hence, we have
\begin{align}
A_1 \ge & \left( \sum_{j\in \mathbb{N}_k} 2^{-jq}\|\triangle_j  \Delta^{-1} ( e^{t\Delta }-1) [\partial_1 (u^0_1 u^0_1) + \partial_2 (u^0_1 u^0_2)] \|^q_\infty \right)^{1/q}  \nonumber\\
& - \left( \sum_{j\in \mathbb{N}_k} 2^{-jq}\|\triangle_j   [ \widetilde{R}_1(u^0_1,  u^0_1) + \widetilde{R}_2 (u^0_1, u^0_2)] \|^q_\infty \right)^{1/q}  := A_{11} -A_{12}. \label{2.21}
\end{align}
In view of Taylor's expansion $e^x=\sum_{r\ge 0} x^r/r! $, we have
\begin{align}
A_{11} \ge & t \left( \sum_{j\in \mathbb{N}_k} 2^{-jq}\|\triangle_j   [\partial_1 (u^0_1 u^0_1) + \partial_2 (u^0_1 u^0_2)] \|^q_\infty \right)^{1/q}  \nonumber\\
& -  t \sum_{r\ge 2}  \frac{1}{r!} \left( \sum_{j\in \mathbb{N}_k} 2^{-jq}\| \triangle_j (t\Delta)^{r-1}  [\partial_1 (u^0_1 u^0_1) + \partial_2 (u^0_1 u^0_2)] \|^q_\infty \right)^{1/q} \nonumber\\
\ge & t \left( \sum_{j\in \mathbb{N}_k} 2^{-jq}\|\triangle_j   (\partial_1- \partial_2)(u^0_1 u^0_1)   \|^q_\infty \right)^{1/q}  \label{2.22} \\
 & - t \left( \sum_{j\in \mathbb{N}_k} 2^{-jq}\|\triangle_j     \partial_2 (u^0_1+ u^0_2) u^0_1   \|^q_\infty \right)^{1/q}   \label{2.23}\\
& -  t \sum_{r\ge 2} \frac{1}{r!} \left( \sum_{j\in \mathbb{N}_k} 2^{-jq}\|\triangle_j (t\Delta)^{r-1} [\partial_1 (u^0_1 u^0_1) + \partial_2 (u^0_1 u^0_2)] \|^q_\infty \right)^{1/q}. \label{2.24}
\end{align}
We will show that $A_{12}$, $\eqref{2.23}$ and $\eqref{2.24}$ are much less than  $\eqref{2.22}$. First, we have

\begin{lem}\label{NSlem2}
Let   $1\leq q\le \infty$,  $t=\eta 2^{-2k}$.  Then
\begin{align}
 t \left( \sum_{j\in \mathbb{N}_k} 2^{-jq}\|\triangle_j   (\partial_1- \partial_2)(u^0_1 u^0_1)   \|^q_\infty \right)^{1/q} \gtrsim \eta \varepsilon k^{1/q}.  \label{2.22bound}
\end{align}
\end{lem}
{\bf Proof.} For the sake of convenience, we denote
\begin{align}
     & \widehat{\Phi^{++}_l}   =  e^{{\rm i}\xi a_l} \varrho (\xi- c_k- b_l) ,  \quad  \widehat{\Phi^{+-}_l}   =  e^{{\rm i}\xi a_l} \varrho (\xi- c_k+ b_l)   \label{NSlem2-2} \\
& \widehat{\Phi^{-+}_l}   =  e^{{\rm i}\xi a_l} \varrho (\xi+ c_k- b_l) ,  \quad  \widehat{\Phi^{--}_l}   =  e^{{\rm i}\xi a_l} \varrho (\xi + c_k+ b_l). \label{NSlem2-3}
\end{align}
Let us observe that
\begin{align}
   \widehat{u^0_1} * \widehat{u^0_1} = &   2^{2k} \sum_{l,m \in \mathbb{N}_k} (\widehat{\Phi^{++}_l} * \widehat{\Phi^{++}_m} + \widehat{\Phi^{+-}_l} * \widehat{\Phi^{+-}_m}) \nonumber\\
   & +  2^{2k} \sum_{l,m \in \mathbb{N}_k} (\widehat{\Phi^{-+}_l} * \widehat{\Phi^{-+}_m} + \widehat{\Phi^{--}_l} * \widehat{\Phi^{--}_m}) \nonumber\\
& +  2^{2k+1} \sum_{l,m \in \mathbb{N}_k} (\widehat{\Phi^{++}_l} * \widehat{\Phi^{+-}_m} + \widehat{\Phi^{-+}_l} * \widehat{\Phi^{--}_m}) \nonumber\\
& +  2^{2k+1} \sum_{l,m \in \mathbb{N}_k} (\widehat{\Phi^{++}_l} * \widehat{\Phi^{--}_m} + \widehat{\Phi^{+-}_l} * \widehat{\Phi^{-+}_m}) \nonumber\\
 & +  2^{2k+1} \sum_{l,m \in \mathbb{N}_k} (\widehat{\Phi^{++}_l} * \widehat{\Phi^{-+}_m} + \widehat{\Phi^{--}_l} * \widehat{\Phi^{+-}_m}) \nonumber\\
 := &   \widehat{U}_1+...+\widehat{U}_5.  \nonumber\\\label{NSlem2-1}
\end{align}
Since ${\rm supp} \varrho (\cdot -a) * \varrho (\cdot -b) \subset B(a+b, 1)$, we see that
\begin{align}
     & {\rm supp} \widehat{U_1}  \cup
        {\rm supp}  \widehat{U_2} \cup
       {\rm supp}  \widehat{U_3} \subset B(2 c_k, \ 2^{2+k/2}) \cup B(-2 c_k, \ 2^{2+k/2}).    \label{NSlem2-5}
\end{align}
It follows that
$
     \triangle_j     (U_1+U_2+U_3 )=0, \ j\in \mathbb{N}_k.
$
Moreover, noticing that ${\rm supp} \widehat{\Phi^{++}_l} * \widehat{\Phi^{--}_m} \subset B(b_l -b_m, \ 1) $, we see that ${\rm supp} \widehat{\Phi^{++}_l} * \widehat{\Phi^{--}_l} \subset B(0,1) $ and  ${\rm supp} \widehat{\Phi^{++}_l} * \widehat{\Phi^{--}_m} \subset  \{\xi: \ 2^{l-2} \le |\xi| < 2^{l-1}\} $ if $m<l$.  Hence, we have
$
     \triangle_j      U_4=0, \ j\in \mathbb{N}_k.
$
  It follows that for any $j\in \mathbb{N}_k$,
\begin{align}
     \triangle_j (u^0_1 u^0_1) =   \triangle_j   U_5 .   \label{NSlem2.29}
\end{align}
Let us rewrite   $U_5$ as
\begin{align}
\widehat{ U_5}
=  & 2^{2k+1} \sum_{l  \in \mathbb{N}_k} (\widehat{\Phi^{++}_l} * \widehat{\Phi^{-+}_l} + \widehat{\Phi^{+-}_l} * \widehat{\Phi^{--}_l}) \nonumber\\
  & + 2^{2k+1} \sum_{l,m \in \mathbb{N}_k, l\neq m} (\widehat{\Phi^{++}_l} * \widehat{\Phi^{-+}_m} + \widehat{\Phi^{+-}_l} * \widehat{\Phi^{--}_m}):= \widehat{U_{51}} +\widehat{U_{52}}. \label{NSlem2.31}
\end{align}
 We easily see that $\triangle_j( \Phi^{++}_l   \Phi^{-+}_l)=0$ if $l \neq j$. Moreover, noticing that $\varphi_j(\xi) =1$ for $|\xi|\in [3\cdot 2^{j-2}, 5\cdot 2^{j-2}]$, we have $\triangle_j ( \Phi^{++}_j \Phi^{-+}_j )= \Phi^{++}_j \Phi^{-+}_j$.   So,
\begin{align}
     \triangle_j U_5 =  2^{2k+1} (\Phi^{++}_j    \Phi^{-+}_j + \Phi^{+-}_j    \Phi^{--}_j)     +  \triangle_j (U_{52}) .   \label{2.30}
\end{align}
In view of \eqref{NSlem2.29} and \eqref{2.30},
\begin{align}
 & t\left( \sum_{j\in \mathbb{N}_k} 2^{-jq}\|\triangle_j   (\partial_1- \partial_2)(u^0_1 u^0_1)   \|^q_\infty  \right)^{1/q}  \nonumber\\
 & \ \  \geqslant  \eta \left( \sum_{j\in \mathbb{N}_k} 2^{-jq} \left\|(\partial_1- \partial_2)   (\Phi^{++}_j    \Phi^{-+}_j + \Phi^{+-}_j    \Phi^{--}_j)   \right \|^q_\infty  \right)^{1/q} \nonumber\\
  & \ \  \ \ -  \eta \left( \sum_{j\in \mathbb{N}_k} 2^{-jq} \left\|(\partial_1- \partial_2) \triangle_j   U_{52}  \right \|^q_\infty  \right)^{1/q}.
    \label{NSlem2-6}
\end{align}
Noticing that
\begin{align}
 \Phi^{+\pm}_l    =   e^{{\rm i} (x + a_l)(c_k \pm b_l)}      \check{\varrho}  (x+a_l), \quad  \Phi^{- \pm}_l    =   e^{{\rm i} (x + a_l)(-c_k \pm b_l)}      \check{\varrho}  (x+a_l),
 \end{align}
one sees that
\begin{align}
\Phi^{++}_j    \Phi^{-+}_j + \Phi^{+-}_j    \Phi^{--}_j      =   \cos {[ a_j (x + a_j)] }   \  ( \check{\varrho}  (x+a_j))^2.
 \end{align}
It follows that
 \begin{align}
(\partial_1 -\partial_2) [\Phi^{++}_j    \Phi^{-+}_j + \Phi^{+-}_j    \Phi^{--}_j  ]    = & - \varepsilon 2^j  \sin {[ a_j (x + a_j)] }   \  ( \check{\varrho}  (x+a_j))^2 \nonumber\\
 & +    \cos {[ a_j (x + a_j)] }   \ (\partial_1 -\partial_2) ( \check{\varrho}  (x+a_j))^2.
 \end{align}
By choosing $x= -a_j + \pi/2a_j$, we have
\begin{align}
\|(\partial_1 -\partial_2) (\Phi^{++}_j    \Phi^{-+}_j + \Phi^{+-}_j    \Phi^{--}_j ) \|_\infty   \ge     \varepsilon 2^j    \    \check{\varrho}  (\pi/2a_j)^2 \ge  \varepsilon 2^{j-1}        \check{\varrho}  (0)^2, \label{lem2.37}
 \end{align}
where we have applied the continuity of $\check{\varrho}$ and assume that  $\check{\varrho}  (\pi/2a_j)^2 \ge  \check{\varrho}  (0)^2/2$.
By \eqref{lem2.37},
\begin{align}
  \eta \left( \sum_{j\in \mathbb{N}_k} 2^{-jq} \left\|(\partial_1- \partial_2)   (\Phi^{++}_j    \Phi^{-+}_j + \Phi^{+-}_j    \Phi^{--}_j)   \right \|^q_\infty  \right)^{1/q} \gtrsim  \eta \varepsilon k^{1/q}.
\label{NSlem2-7}
\end{align}
For any $j, l, m \in \mathbb{N}_k $ with $l<m$, if $m\neq j$, then we have $\triangle_j (\Phi^{++}_m \Phi^{-+}_l)=0$. Hence, applying the multiplier estimates, we have
\begin{align}
  \left\| \partial_i \triangle_j \left(\sum_{m,l\in \mathbb{N}_k; m> l} \Phi^{++}_m \Phi^{-+}_l \right)  \right \|_\infty  \lesssim 2^j \left\|  \sum_{l<j,  l\in \mathbb{N}_k} \Phi^{++}_j \Phi^{-+}_l    \right \|_\infty, \ \ i=1,2.
    \label{NSlem2-9}
\end{align}
Using the rapid decay \eqref{rapiddecay}, we have
\begin{align}
   \left\|  \sum_{l<j,  l\in \mathbb{N}_k} \Phi^{++}_j \Phi^{-+}_l    \right \|_\infty & \lesssim   \left\|  \sum_{l<j,  l\in \mathbb{N}_k}  (1+|x+a_j|)^{-N} (1+|x+a_l|)^{-N}    \right \|_\infty \nonumber\\
& \lesssim   \left\|  \sum_{l<j,  l\in \mathbb{N}_k}  (1+|x|)^{-N} (1+|x+a_l-a_j|)^{-N}    \right \|_\infty .
    \label{NSlem2-10}
\end{align}
In \eqref{NSlem2-10}, by separating $\mathbb{R}^n$ into two different regions $\{x: |x|\le 3\cdot 2^{j-3}\}$ and $\{x: |x|> 3\cdot 2^{j-3}\}$, we easily see that the right hand side of \eqref{NSlem2-10} can be bounded by $j 2^{-Nj}$. So,
\begin{align}
  \left\| \partial_i \triangle_j \left(\sum_{m,l\in \mathbb{N}_k; m> l} \Phi^{++}_m \Phi^{-+}_l \right)  \right \|_\infty  \lesssim 2^j j 2^{-Nj}, \ \ i=1,2 .
    \label{NSlem2-11}
\end{align}
By symmetry, \eqref{NSlem2-11} also holds if one substitutes  the summation $\sum_{m,l\in \mathbb{N}_k; m> l}$ by  $\sum_{m,l\in \mathbb{N}_k; m< l}$.   It follows from \eqref{NSlem2-11} that
\begin{align}
 \left( \sum_{j\in \mathbb{N}_k} 2^{-jq} \left\| \partial_i \triangle_j \left(\sum_{m,l\in \mathbb{N}_k; m\neq l} \Phi^{++}_m \Phi^{-+}_l \right)  \right \|^q_\infty  \right)^{1/q} \lesssim k 2^{-k}, \ \ i=1,2.
    \label{NSlem2-12}
\end{align}
Using the same way as in \eqref{NSlem2-12}, we can estimate another term in $U_{52}$ and
\begin{align}
 \left( \sum_{j\in \mathbb{N}_k} 2^{-jq} \left\| \partial_i \triangle_j U_{52}  \right \|^q_\infty  \right)^{1/q} \lesssim k 2^{-k}, \ \ i=1,2.
    \label{NSlem2-13}
\end{align}

Collecting the estimates as in \eqref{NSlem2-6},  \eqref{NSlem2-7}  and \eqref{NSlem2-13}, we have
\begin{align}
 t \left( \sum_{j\in \mathbb{N}_k} 2^{-jq}\|\triangle_j   (\partial_1- \partial_2)(u^0_1 u^0_1)   \|^q_\infty \right)^{1/q} \gtrsim \eta \varepsilon k^{1/q}   - C \eta 2^{-k/4}  - C\eta k 2^{-k}.  \label{2.22boundin}
\end{align}
By choosing $k\gg -\log \varepsilon$, we have the lower bound as desired in \eqref{2.22bound}. $\hfill\Box$

\begin{lem}\label{NSlem3}
Let   $1\leq q\le \infty$,  $t=\eta 2^{-2k}$.  Then
\begin{align}
 t \left( \sum_{j\in \mathbb{N}_k} 2^{-jq}\|\triangle_j     \partial_2 (u^0_1+ u^0_2) u^0_1   \|^q_\infty \right)^{1/q}  \lesssim \eta\varepsilon k^{1+1/q} 2^{-k/2}     \label{2.39}
\end{align}
\end{lem}
{\bf Proof. } By \eqref{u0}, we have
\begin{align}
u^0_1 + u^0_2 &= 2^k  \sum_{l\in \mathbb{N}_k}  \mathscr{F}^{-1} \frac{\xi_2-\xi_1}{\xi_2}   (\widehat{\Phi^+_l}(\xi)+ \widehat{\Phi^-_l}(\xi)). \label{2.40}
  \end{align}
By the multiplier estimate and \eqref{2.11},
\begin{align}
  & t \left( \sum_{j\in \mathbb{N}_k} 2^{-jq}\|\triangle_j     \partial_2 (u^0_1+ u^0_2) u^0_1   \|^q_\infty \right)^{1/q}  \nonumber\\
  & \quad
   \lesssim   t k^{1/q}   \| u^0_1+ u^0_2 \|_\infty \| u^0_1   \|_\infty  \lesssim   \eta 2^{-k} k^{1/q}   \| u^0_1+ u^0_2 \|_\infty.  \label{2.41}
\end{align}
It follows from \eqref{2.40},
\eqref{2.13} and \eqref{2.14} that
\begin{align}
   \| u^0_1+ u^0_2 \|_\infty  & \lesssim  2^{k}  \sum_{l\in \mathbb{N}_k} \sum_{\lambda,\mu=\pm 1} \left\| \mathscr{F}^{-1} \left(\frac{\xi_2-\xi_1}{\xi_2} \psi(\xi + \lambda c_k - \mu b_l)\right)\right\|_1  \nonumber\\
   & \lesssim \varepsilon k 2^{k/2}.     \label{2.42}
\end{align}
In view of \eqref{2.41} and \eqref{2.42}, we have \eqref{2.39}. $\hfill\Box$

\begin{lem}\label{NSlem4}
Let   $1\leq q\le \infty$,  $t=\eta 2^{-2k}$.  Then
\begin{align}
  t \sum_{r\ge 2} \frac{1}{r!} \left( \sum_{j\in \mathbb{N}_k} 2^{-jq}\|   (t\Delta)^{r-1} \triangle_j     \partial_i (u^0_\alpha u^0_\beta)  \|^q_\infty \right)^{1/q} \lesssim \eta^2    k^{1/q} 2^{-k}, \ i,\alpha,\beta=1,2. \label{2.43}
\end{align}
\end{lem}
{\bf Proof.} Using Bernstein's multiplier estimates, we have for $i=1,2$,
\begin{align}
 \|\partial_i \mathscr{F}^{-1} (t|\xi|^2)^{r-1} \varphi_j  \mathscr{F} f \|_\infty & \lesssim 2^j (t2^{2j})^{r-1}  \| \mathscr{F}^{-1} |\xi|^{2(r-1)} \varphi \|_1 \|f\|_\infty \nonumber\\
  & \lesssim 2^j (t2^{2j})^{r-1} 4^r r^n    \|f\|_\infty.
 \label{2.43}
\end{align}
Noticing that $j\in \mathbb{N}_k$ implies that $j\le k/2$,  one has that for $i=1,2$,
\begin{align}
&  t \sum_{r\ge 2} \frac{1}{r!} \left( \sum_{j\in \mathbb{N}_k} 2^{-jq}\| \mathscr{F}^{-1} (t|\xi|^2)^{r-1} \varphi_j  \mathscr{F}     \partial_i (u^0_\alpha u^0_\beta)  \|^q_\infty \right)^{1/q} \nonumber\\
& \quad \lesssim  t \sum_{r\ge 2} \frac{ 4^r r^n }{r!} \left( \sum_{j\in \mathbb{N}_k}  ((t2^{2j})^{q(r-1)}   \| u^0_\alpha u^0_\beta   \|^q_\infty \right)^{1/q} \nonumber\\
& \quad \lesssim \eta \sum_{r\ge 2} \frac{ 4^r r^n }{r!} k^{1/q} (t2^{k})^{ (r-1)} \lesssim \eta^2   k^{1/q} 2^{-k}.
 \label{2.44}
\end{align}
Collecting Lemmas \ref{NSlem2}--\ref{NSlem4}, we immediately have

\begin{lem}\label{NSlem5}
Let   $1\leq q\le \infty$,  $t=\eta 2^{-2k}$, $0<\eta \le \varepsilon^2 \ll 1, \ k\gg -\log \varepsilon$.  Then
\begin{align}
  A_{11}  \gtrsim \eta \varepsilon    k^{1/q}.  \label{2.45}
\end{align}
\end{lem}
In the following we will show that $A_{12}$ is much less than $A_{11}$. One can rewrite $\widetilde{R}_i$ as
\begin{align}
\widetilde{R}_i(f,g) = \int^t_0  e^{(t-\tau)\Delta }  \partial_i  [e^{\tau\Delta } f \ ( e^{\tau\Delta }-1) g + g( e^{\tau\Delta }-1)f ] d\tau. \label{2.46}
\end{align}
\begin{lem}\label{NSlem6}
Let   $1\leq q\le \infty$,  $t=\eta 2^{-2k}$.  Then
\begin{align}
    \left( \sum_{j\in \mathbb{N}_k} 2^{-jq}\| \triangle_j \widetilde{R}_i(u^0_\alpha, u^0_\beta)\|^q_\infty \right)^{1/q} \lesssim \eta^2    k^{1/q}, \ \ i, \alpha, \beta=1,2. \label{2.47}
\end{align}
\end{lem}
{\bf Proof.} Using the fact that $e^{t\Delta}: L^\infty \to L^\infty$, we have
\begin{align}
 &   \left( \sum_{j\in \mathbb{N}_k} 2^{-jq}\| \triangle_j \widetilde{R}_1(u^0_1, u^0_1)  \|^q_\infty \right)^{1/q}  \nonumber\\
 & \quad \lesssim t \sup_{\tau\in [0,t]} \left( \sum_{j\in \mathbb{N}_k} \| \triangle_j [e^{\tau\Delta } u^0_1 \ ( e^{\tau\Delta }-1) u^0_1 + u^0_1( e^{\tau\Delta }-1)u^0_1 ]  \|^q_\infty \right)^{1/q} \nonumber\\
  & \quad \lesssim t  k^{1/q} \sup_{\tau\in [0,t]}   \|     u^0_1 \|_\infty \| ( e^{\tau\Delta }-1) u^0_1   \|_\infty  \nonumber\\
   & \quad \lesssim \eta  k^{1/q} 2^{-k} \sup_{\tau\in [0,t]}    \| ( e^{\tau\Delta }-1) u^0_1   \|_\infty .
  \label{2.49}
\end{align}
Using Taylor's expansion, one has that
\begin{align}
\| ( e^{\tau\Delta }-1) u^0_1   \|_\infty  \leqslant \sum^\infty_{r=1}   \frac{\tau^r}{r!}\|\mathscr{F}^{-1}|\xi|^{2r} \mathscr{F} u^0_1  \|_\infty. \label{2.50}
\end{align}
Since ${\rm supp} \ \widehat{u^0_1} \subset \{\xi: 2^k \le |\xi| \le 2^k+ 2^{k/2}+1\}$, we see that
\begin{align}
 \|\mathscr{F}^{-1}|\xi|^{2r} \mathscr{F} u^0_1  \|_\infty  & \leqslant \|\mathscr{F}^{-1}(|\xi|^{2r}\varphi_k)\|_1 \|u^0_1\|_\infty \nonumber\\
 & \leqslant 2^{2kr}\|\mathscr{F}^{-1}(|\xi|^{2r}\varphi)\|_1 \|u^0_1\|_\infty \lesssim r^n 4^r 2^{2kr+1} \lesssim   8^r 2^{2kr+1}.
 \label{2.51}
\end{align}
It follows from \eqref{2.50} and \eqref{2.51} that
\begin{align}
\| ( e^{\tau\Delta }-1) u^0_1   \|_\infty  \leqslant  C2^k (e^{8\tau 2^{2k}}-1)  \leqslant  C2^k (e^{8\eta}-1) \leqslant C\eta 2^k.  \label{2.52}
\end{align}
By \eqref{2.49} and \eqref{2.52}, we have the result, as desired. The other cases can be proven in a similar way.  $\hfill \Box$

If we take $\eta=\varepsilon^2,$ then we see that $\eta^2 k^{1/q} \ll \eta \varepsilon k^{1/q}$. So, $A_{11} \gg A_{12}$. In the following we need to control $A_2$ and show that it is much less than $A_1$. For convenience, we use the same notations as in \eqref{2.18} and \eqref{2.19}.  One can rewrite $F_2$ as
\begin{align}
F_2(t,x) = \sum^2_{\alpha,\beta=1} \frac{\partial_\alpha\partial_\beta}{\Delta} R_1(u^0_\alpha, u^0_\beta). \label{2.53}
\end{align}

\begin{lem}\label{NSlem7}
Let   $1\leq q\le \infty$, $0<\varepsilon \ll 1,$  $t=\eta 2^{-2k}$, $\eta \le \varepsilon ^2$.  Then
\begin{align}
    A_2  \lesssim \eta \varepsilon^2  k^{1/q}. \label{2.54}
\end{align}
\end{lem}
{Proof.} Recalling that
\begin{align}
    R_1(f,g) = (e^{t\Delta}-1) \Delta^{-1}\partial_1 (fg) + \widetilde{R}_1(f,g), \label{2.55}
\end{align}
we see that
\begin{align}
   A_2 \le & \sum^2_{\alpha, \beta=1} \left( \sum_{j\in \mathbb{N}_k} 2^{-jq} \left\|\triangle_j   \frac{\partial_\alpha\partial_\beta}{\Delta} (e^{t\Delta}-1) \Delta^{-1} \partial_1 (u^0_\alpha  u^0_\beta) \right\|^q_\infty \right)^{1/q} \nonumber\\
   & + \sum^2_{\alpha, \beta=1} \left( \sum_{j\in \mathbb{N}_k} 2^{-jq} \left\|\triangle_j   \frac{\partial_\alpha\partial_\beta}{\Delta} \widetilde{R}_1(u^0_\alpha, u^0_\beta)  \right\|^q_\infty \right)^{1/q} \nonumber\\
   := & A_{21} + A_{22}. \label{2.56}
\end{align}
The estimate of $A_{22} $ is very easy. Noticing that
\begin{align}
     \left\|\triangle_j   \frac{\partial_\alpha\partial_\beta}{\Delta} f \right\|_\infty
  \lesssim \sum^{1}_{\ell =-1} \left\|\mathscr{F}^{-1} \frac{\xi_\alpha\xi_\beta}{|\xi|^2} \varphi_{j+\ell}(\xi) \right\|_1 \|\triangle_j f\|_\infty \lesssim \|\triangle_j f\|_\infty,   \label{2.57}
\end{align}
from Lemma \ref{NSlem6} we immediately have
\begin{align}
  A_{22}  \lesssim  \sum^2_{\alpha, \beta=1} \left( \sum_{j\in \mathbb{N}_k} 2^{-jq} \left\|\triangle_j   \widetilde{R}_1(u^0_\alpha, u^0_\beta)  \right\|^q_\infty \right)^{1/q} \lesssim \eta^2 k^{1/q}. \label{2.58}
\end{align}
Now we estimate $A_{21}$. Using Taylor's expansion and \eqref{2.57},
\begin{align}
   \left\|\triangle_j \frac{\partial_\alpha\partial_\beta}{\Delta} (e^{t\Delta}-1) \Delta^{-1} \partial_1(u^0_\alpha  u^0_\beta) \right\|_\infty
  & \lesssim t \left\|\triangle_j \frac{\partial_\alpha\partial_\beta}{\Delta} \partial_1 (u^0_\alpha  u^0_\beta) \right\|_\infty \nonumber\\
  & \quad + t \sum^\infty_{r=2} \frac{1}{r!}  \left\|\triangle_j (t\Delta)^{r-1} \partial_1 (u^0_\alpha  u^0_\beta) \right\|_\infty.
   \label{2.59}
\end{align}
Hence, we have
\begin{align}
   A_{21} \le & t \sum^2_{\alpha, \beta=1} \left( \sum_{j\in \mathbb{N}_k} 2^{-jq} \left\|\triangle_j   \frac{\partial_\alpha\partial_\beta}{\Delta} \partial_1  (u^0_\alpha  u^0_\beta) \right\|^q_\infty \right)^{1/q} \nonumber\\
   & + t \sum^2_{\alpha, \beta=1} \sum^\infty_{r=2} \frac{1}{r!}\left( \sum_{j\in \mathbb{N}_k} 2^{-jq} \left\|\triangle_j   (t\Delta)^{r-1} \partial_1 (u^0_\alpha  u^0_\beta)  \right\|^q_\infty \right)^{1/q} \nonumber\\
   := & \sum_{\alpha, \beta=1,2} A^{\alpha\beta}_{211} + A_{212}. \label{2.60}
\end{align}
By Lemma \ref{NSlem4},
\begin{align}
   A_{212} \lesssim \eta^2 k^{1/q} 2^{-k}. \label{2.61}
\end{align}
So, it suffices to bound $A^{\alpha\beta}_{211}$. We divide the estimates of $A^{\alpha\beta}_{211}$ into the following three cases.

{\it Case} 1. $\alpha=\beta =1$. From \eqref{NSlem2.29}, \eqref{NSlem2.31} and \eqref{2.30}, it follows that for $j\in \mathbb{N}_k$,
\begin{align}
  \triangle_j    (u^0_1 u^0_1)   =  2^{2k+1} (\Phi^{++}_j \Phi^{-+}_j + \Phi^{+-}_j    \Phi^{--}_j)     +  \triangle_j U_{52},
    \label{2.62}
\end{align}
we see that
\begin{align}
   A^{11}_{211} = & t   \left( \sum_{j\in \mathbb{N}_k} 2^{-jq} \left\|\triangle_j   \frac{\partial^3_1 }{\Delta}    (u^0_1  u^0_1) \right\|^q_\infty \right)^{1/q} \nonumber\\
  \le  & 2\eta   \left( \sum_{j\in \mathbb{N}_k} 2^{-jq} \left\|   \frac{\partial^3_1 }{\Delta}   \triangle_j \left(\Phi^{++}_j \Phi^{-+}_j + \Phi^{+-}_j    \Phi^{--}_j\right) \right\|^q_\infty \right)^{1/q} \nonumber\\
  & + 2\eta   \left( \sum_{j\in \mathbb{N}_k} 2^{-jq} \left\|   \frac{\partial^3_1 }{\Delta}   \triangle_j  U_{52}   \right\|^q_\infty \right)^{1/q} \nonumber\\
  := & B_1+B_2. \label{2.63}
\end{align}
Using the Bernstein's multiplier estimate and   \eqref{NSlem2-13},
\begin{align}
  B_2 \lesssim  \eta   \left( \sum_{j\in \mathbb{N}_k} 2^{-jq} \left\|    \triangle_j \partial_1    U_{52}   \right\|^q_\infty \right)^{1/q} \lesssim \eta k 2^{-k}.   \label{2.64}
\end{align}
For the estimate of $B_1$, noticing that ${\rm supp} \ \widehat{\Phi^{++}_j} * \widehat{\Phi^{-+}_j} \subset B(a_j,1)$, we have from the multiplier estimate,
\begin{align}
  \left\|  \triangle_j \frac{\partial^3_1 }{\Delta}   \left(  \Phi^{++}_j \Phi^{-+}_j \right) \right\|_\infty & \lesssim  2^j \left\| \mathscr{F}^{-1} \frac{\xi^2_1 }{|\xi|^2} \psi(\xi-a_j) \right\|_1  \left\|  \Phi^{++}_j \Phi^{-+}_j   \right\|_\infty  \nonumber\\
  & \lesssim   2^j \left\| \mathscr{F}^{-1} \frac{(\xi_1 +\varepsilon 2^j)^2 }{|\xi+a_j|^2} \psi(\xi) \right\|_1.
  \label{2.65}
\end{align}
Since the third coordinate of $a_j$ is larger than $2^j/2$, we have from the Bernstein's multiplier estimate that ($j\ge k/4\gg -\log \varepsilon$)
\begin{align}
    \left\| \mathscr{F}^{-1} \frac{(\xi_1 +\varepsilon 2^j)^2 }{|\xi+a_j|^2} \psi(\xi) \right\|_1 \lesssim \varepsilon^2.
  \label{2.66}
\end{align}
Therefore,
\begin{align}
    B_1  \lesssim \eta \varepsilon^2 k^{1/q}.
  \label{2.67}
\end{align}
Summarizing \eqref{2.64} and \eqref{2.67}, one has that
\begin{align}
A^{11}_{211}    \lesssim \eta \varepsilon^2 k^{1/q}.
  \label{2.68}
\end{align}

{\it Case} 2. $\alpha=\beta =2$. We can rewrite $u^0_2$ as
\begin{align}
 u^0_2   &    = -u^0_1  +  \widetilde{u^0_2}, \ \  \widetilde{u^0_2}=  2^k \sum_{l\in \mathbb{N}_k} \mathscr{F}^{-1} \frac{\xi_2-\xi_1}{\xi_2} ( \widehat{\Phi^+_l}(\xi) + \widehat{\Phi^-_l}(\xi)).
   \label{2.69}
  \end{align}
It follows that
\begin{align}
 u^0_2  u^0_2 &    =  u^0_1 u^0_1  +  \widetilde{u^0_2} \widetilde{u^0_2} - 2 u^0_1\widetilde{u^0_2}  \label{2.70}
  \end{align}
and
\begin{align}
   A^{22}_{211} \le  & t   \left( \sum_{j\in \mathbb{N}_k} 2^{-jq} \left\|\triangle_j   \frac{\partial_1 \partial^2_2 }{\Delta}    (u^0_1  u^0_1) \right\|^q_\infty \right)^{1/q} \nonumber\\
     &  + 2 t   \left( \sum_{j\in \mathbb{N}_k} 2^{-jq} \left\|\triangle_j   \frac{\partial_1 \partial^2_2 }{\Delta}    (u^0_1 \widetilde{u^0_2}  ) \right\|^q_\infty \right)^{1/q} \nonumber\\
     & +  t   \left( \sum_{j\in \mathbb{N}_k} 2^{-jq} \left\|\triangle_j   \frac{\partial_1 \partial^2_2}{\Delta}    (\widetilde{u^0_2} \widetilde{u^0_2}  ) \right\|^q_\infty \right)^{1/q}:= I+II+III.  \label{2.71}
\end{align}
Applying the same way as in the estimate of $B_1$, one has that
\begin{align}
I  \lesssim \eta \varepsilon^2 k^{1/q}.
  \label{2.72}
\end{align}
By Bernstein's multiplier estimate and Lemma \ref{NSlem3}, we have
\begin{align}
  II \lesssim  t   \left( \sum_{j\in \mathbb{N}_k} 2^{-jq} \left\|\triangle_j    \partial_1 (u^0_1 \widetilde{u^0_2}  ) \right\|^q_\infty \right)^{1/q}  \lesssim \eta\varepsilon k^{1+1/q} 2^{-k/2}. \label{2.73}
\end{align}
The estimate of $III$ is easier than that of $II$ and we have
\begin{align}
  III \lesssim  \eta\varepsilon k^{2+1/q} 2^{-k}. \label{2.74}
\end{align}
Collecting the estimates of \eqref{2.72},  \eqref{2.73} and \eqref{2.74}, we have
\begin{align}
   A^{22}_{211} \le  \eta\varepsilon^2 k^{1/q} .  \label{2.75}
\end{align}

{\it Case} 3. $\alpha=1, \ \beta =2$. Similar to Case 2, we can show that $ A^{12}_{211}$ has the same bound as $ A^{22}_{211}$ and we omit the details of the proof.

Up to now, we have shown that
\begin{align}
  \sum_{\alpha,\beta=1,2} A^{\alpha\beta}_{211} \le  \eta\varepsilon^2 k^{1/q} .  \label{2.76}
\end{align}
By \eqref{2.60}, \eqref{2.61} and \eqref{2.76}, we have
\begin{align}
  A_{21} \le  \eta\varepsilon^2 k^{1/q} .  \label{2.77}
\end{align}
In view of \eqref{2.56}, \eqref{2.58} and \eqref{2.77}, we immediately have the result, as desired. $\hfill\Box$\\

\begin{lem}\label{NSlem8}
Let   $1\leq q\le \infty$,  $t= \varepsilon^2 2^{-2k}$.  Then
\begin{align}
  \left( \sum_{j\in \mathbb{N}_k} 2^{-jq} \left\|\triangle_j  \frac{\partial^2 u_1}{\partial \delta^2}(0,t) \right \|^q_\infty \right)^{1/q} \gtrsim   \varepsilon^3 k^{1/q}.  \label{2.78}
\end{align}
\end{lem}
{\bf Proof.} In view of Lemma \ref{NSlem6} we have $A_{12} \lesssim \varepsilon^4 k^{1/q}$. Using \eqref{2.21} and Lemma \ref{NSlem5}, we get  $A_{1} \ge (c \varepsilon^3-C \varepsilon^4) k^{1/q}$. By Lemma \ref{NSlem7} and \eqref{secondder}, one has that for $t=\varepsilon^2 2^{-2k}$,
\begin{align}
 \left( \sum_{j\in \mathbb{N}_k} 2^{-jq} \left\|\triangle_j  \frac{\partial^2 u_1}{\partial \delta^2}(0,t) \right \|^q_\infty \right)^{1/q}  \ge (c \varepsilon^3-2 C \varepsilon^4) k^{1/q}.  \label{2.78}
\end{align}
We can choose $\varepsilon $ satisfying $100C\varepsilon \le c/2$, the result follows. $\hfill \Box$

\subsection{Estimates on higher order terms of $\delta$} \label{sect2.2}

Using Koch and Tataru's global wellposedness result, we can show that the higher order terms on $\delta$ in Taylor's expansion \eqref{1.8} in $\dot B^{-1}_{\infty,q}$  are much smaller than $\delta^2\varepsilon^3 k^{1/q}$. This fact, together with Lemma \ref{NSlem8} imply the inflation phenomena of the solution of NS in $\dot B^{-1}_{\infty,q}$.   By induction we have for the solution $u(\delta,t)$ of NS,
\begin{align}
\frac{\partial^r u}{\partial \delta^r}(0,t) = \int^t_0 e^{(t-\tau)\Delta} \mathbb{P} {\rm div} \sum^{r-1}_{m=1} {r \choose m} \left( \frac{\partial^m u}{\partial \delta^m}(0,\tau) \otimes \frac{\partial^{r-m} u}{\partial \delta^{r-m}}(0,\tau)\right) d\tau  .  \label{2.80}
\end{align}
Let $X$ be defined in Theorem \ref{NSWellposed}. Let us recall that Koch and Tataru \cite{KoTa01}, Iwabuchi and Nakamura  \cite{Iw13}  obtained the following estimates.
\begin{align}
\left\|\int^t_0 e^{(t-\tau)\Delta} \mathbb{P} {\rm div}   \left( u(\tau) \otimes v(\tau)\right) d\tau \right\|_{Y}  \lesssim \|u\|_X\|v\|_X,  \label{Koch2.80}
\end{align}
where $Y:= X\cap L^\infty(0,\infty; BMO^{-1})$ Applying the integral equation, \eqref{Koch2.80} and Theorem \ref{NSWellposed}, we see that
\begin{align}
\|u (\delta, t)-\delta e^{t\Delta}u_0\|_Y \le  \left\|\int^t_0 e^{(t-\tau)\Delta} \mathbb{P} {\rm div}   \left( u(\tau) \otimes u(\tau)\right) d\tau \right\|_Y  \lesssim \|u\|^2_X \lesssim\delta^2. \label{Koch2.81}
\end{align}
One sees that
\begin{align}
\label{Koch2.82}
  u (\delta, t)- \delta e^{t\Delta}u_0 -  \frac{\delta^2}{2} \frac{\partial^2 u}{\partial \delta^2}(0,t)&  = \int^t_0 e^{(t-\tau)\Delta} \mathbb{P} {\rm div}   \left( u(\tau) \otimes u(\tau) - \delta^2e^{\tau\Delta}u_0 \otimes e^{\tau\Delta}u_0 \right) d\tau.
\end{align}
It follows from \eqref{Koch2.80}, \eqref{Koch2.81}, \eqref{Koch2.82} and Theorem \ref{NSWellposed} that
\begin{align}
\label{Koch2.83}
  \left\|u (\delta, t)- \delta e^{t\Delta}u_0 -  \frac{\delta^2}{2} \frac{\partial^2 u}{\partial \delta^2}(0,t) \right\|_Y \lesssim  \|u\|_X \left\|u (\delta, t)- \delta e^{t\Delta}u_0  \right\|_X \lesssim \delta^3.
\end{align}
Let us fix $t=\varepsilon^2 2^{-2k}$ and denote
\begin{align}
\label{2.81}
\widetilde{u}(\delta, t) = \sum^\infty_{r=3} \frac{\delta^r}{r!} \frac{\partial^r u}{\partial \delta^r}(0,t).
\end{align}
By \eqref{1.8} we see that there exist $\delta_0>0$, such that for any $0<\delta\le \delta_0$,
\begin{align}
\label{2.82}
  u (\delta, t) = \delta e^{t\Delta}u_0  +   \frac{\delta^2}{2} \frac{\partial^2 u}{\partial \delta^2}(0,t)+ \widetilde{u}(\delta, t) \ \ \ \ \mbox{in} \ \ Y.
\end{align}
By  \eqref{Koch2.83}, one sees that for $\delta\le \delta_0$,
\begin{align}
\label{2.82a}
  \|\widetilde{u}(\delta, t) \|_{Y}   \lesssim \delta^3.
\end{align}
This implies that for any $\delta \le \min\{\delta_0, \ \varepsilon^4\}$,
\begin{align}
\label{2.86}
   \left\|  \widetilde{u}(\delta, t) \right\|_{\dot B^{-1}_{\infty,\infty}} \le   C \delta^2 \varepsilon^4, \ \ \forall t>0.
\end{align}
Applying this fact, we can prove Theorem \ref{NSIll}.\\

\subsection{\bf Proof of Theorem \ref{NSIll}: $n\ge 3$.}

Let $0<\varepsilon \ll 1$ be as in Lemma \ref{NSlem8},  $t=\varepsilon^2 2^{-2k}  $ and $\delta \le \min\{\delta_0, \ \varepsilon^4\}$ be the same one as in \eqref{2.86}.   Using   \eqref{2.82}, we have
\begin{align}
\left\|   u_1  (\delta,t) \right \|_{\dot B^{-1}_{\infty,1}}  & \ge
    \left( \sum_{j\in \mathbb{N}_k} 2^{-jq} \left\|\triangle_j    u_1 (\delta,t) \right \|^q_\infty \right)^{1/q}  \nonumber\\
    & \ge \frac{\delta^2}{2}
   \left( \sum_{j\in \mathbb{N}_k} 2^{-jq} \left\|\triangle_j  \frac{\partial^2 u_1}{\partial \delta^2}(0,t) \right \|^q_\infty \right)^{1/q} \nonumber\\
 & \quad -\delta \|e^{t\Delta}u_0\|_{\dot B^{-1}_{\infty,q}} -
   \left( \sum_{j\in \mathbb{N}_k} 2^{-jq} \left\|\triangle_j   \widetilde{u}_1 (\delta,t) \right \|^q_\infty \right)^{1/q}.    \label{2.87}
\end{align}
 Applying \eqref{2.86} and noticing that $\mathbb{N}_k$ has at most $k$ many indices,   we have
\begin{align}
    \left( \sum_{j\in \mathbb{N}_k} 2^{-jq} \left\|\triangle_j   \widetilde{u}_1 (\delta,t) \right \|^q_\infty \right)^{1/q}
 &  \le   k^{1/q} \left\|  \widetilde{u}_1 (\delta,t) \right \|_{\dot B^{-1}_{\infty,\infty}}
    \le C  k^{1/q} \delta^2 \varepsilon^4 .  \label{2.88}
\end{align}
Obviously, we have
\begin{align}
\|e^{t\Delta}u_0\|_{\dot B^{-1}_{\infty,q}} \le C.  \label{2.89}
\end{align}
By Lemma \ref{NSlem8}, \eqref{2.87}, \eqref{2.88} and \eqref{2.89}, we have
\begin{align}
\left\|   u_1  (\delta,t) \right \|_{\dot B^{-1}_{\infty,1}}    \ge c \delta^2  \varepsilon^3  k^{1/q} -C \delta - C \delta^2 \varepsilon^4 k^{1/q}.
        \label{2.90}
\end{align}
Recalling that $\varepsilon$ has been chosen as in \eqref{2.78},  we immediately have for  $k\gg \varepsilon^{-3q} \delta^{-q}$,
\begin{align}
\left\|   u_1  (\delta,t) \right \|_{\dot B^{-1}_{\infty,1}}    \gtrsim     \varepsilon^{3} \delta^2  k^{1/q}.
        \label{2.99}
\end{align}
This finishes the proof of Theorem \ref{NSIll} in the case $n\ge 3$.  $\hfill\Box$

\section{Proof of Theorem \ref{NSIll}, $n=2$} \label{sect3}

Noticing that in Section \ref{sect2} one needs the third coordinates of $a_l, b_l$ are not zero, the proof of higher dimensional cases cannot be straightly applied to the 2D case. However, the nonlinearity in 2D case is simpler than that of higher dimensions. We now sketch the proof in 2D case. We write
\begin{align}
  \partial_t u_j -\Delta u_j  + B_j (u,u)=0, \ \ j=1,2,   \label{2d3.1}
\end{align}
where
\begin{align}
    B_1 (u,u) & : = \sum^2_{i=1} \partial_i (u_i u_1)  -  \partial_1 \Delta^{-1} \sum^2_{l, m=1} \partial_l \partial_m  (u_l u_m) \nonumber\\
    &= \Delta^{-1} (\partial^2_1-\partial^2_2) \partial_2 (u_1u_1) +  \partial_2 (u_1(u_1+u_2)) +  \Delta^{-1}  \partial^2_2  \partial_1 ((u_1)^2 -(u_2)^2) \nonumber\\
    & \quad  -2 \Delta^{-1}  \partial^2_1  \partial_2 (u_1(u_1+u_2)).   \label{2d3.2}
\end{align}
and $B_2(u,u)$ is similar. We have from \eqref{seconditerat} that
\begin{align} \label{2d3.3}
  \frac{\partial^2  u_1}{\partial\delta^2}\Big|_{\delta=0} = \int^t_0  e^{(t-\tau)\Delta }  B_1 (e^{\tau\Delta } u_0, \   e^{\tau\Delta } u_0) d\tau.
\end{align}
Let $0<\varepsilon \ll 1$, $k$ and $\mathbb{N}_k$ be the same ones as in the cases $n\ge 3$. Put
\begin{align}
   c_k =  (\sqrt{2}/2, \ \sqrt{2}/2 )2^k, \ \ a_l = (\varepsilon, \ \sqrt{1-\varepsilon^2}) 2^l, \ \ b_l=a_l/2    \label{2d3.4}
\end{align}
and let $u^0_1$ and $ u^0_2$ be the same ones as \eqref{u0}. Similarly as in \eqref{2.13}, we have
\begin{align}
   \|u^0_1+ u^0_2\|_\infty \lesssim 2^{k/2}.     \label{2d3.5}
\end{align}
Let us  estimate
\begin{align}
  \left(\sum_{j\in \mathbb{N}_k}  2^{-jq} \left\|\triangle_j \frac{\partial^2  u_1}{\partial\delta^2}\Big|_{\delta=0} \right\|^q_\infty \right)^{1/q}.   \label{2d3.6}
\end{align}
Let $t=\eta 2^{-2k}$.  Denote
\begin{align}
\widetilde{R}(t,u_0) =&   \frac{\partial^2  u_1}{\partial\delta^2}\Big|_{\delta=0} - \int^t_0  e^{(t-\tau)\Delta }  B_1 (  u_0, \     u_0) d\tau \nonumber\\
 =&     \int^t_0  e^{(t-\tau)\Delta } ( B_1 (  e^{\tau\Delta} u_0, \   e^{\tau\Delta}  u_0)- B_1 (  u_0, \     u_0)) d\tau.
 \label{3.2d7}
 \end{align}
Similarly as in \eqref{2.46},  every term in $\widetilde{R}(t,u_0)$ contains $(e^{\tau\Delta}-1)$,  using the same way as in Lemma \ref{NSlem6}, we have for $t=\eta 2^{-2k}$,
\begin{align}
   \left(\sum_{j\in \mathbb{N}_k}  2^{-jq} \left\|\triangle_j \widetilde{R}(t,u_0) \right\|^q_\infty \right)^{1/q}  \lesssim \eta^2 k^{1/q}.   \label{2d3.8}
\end{align}
We can write
\begin{align}
 \int^t_0  e^{(t-\tau)\Delta }  B_1 (  u_0, \     u_0) d\tau =&    \Delta^{-1} (e^{t\Delta}-1) B_1(u_0,\ u_0) \nonumber\\
 := & t\sum_{r\ge 1} \frac{(t\Delta)^{r-1}}{r!} \Delta^{-1} (\partial^2_1- \partial^2_2)\partial_2 (u^0_1u^0_1) +  \widetilde{Q}(t, u_0).
 \label{2d3.9}
 \end{align}
Since every term in $\widetilde{Q}(t, u_0)$  contains $u^0_1+u^0_2$, by \eqref{2d3.5} we can repeat the procedures as in Section \ref{sect2} to obtain that
\begin{align}
   \left(\sum_{j\in \mathbb{N}_k}  2^{-jq} \left\|\triangle_j \widetilde{Q}(t,u_0) \right\|^q_\infty \right)^{1/q}  \lesssim \eta k^{1/q} 2^{-k/2}.   \label{2d3.10}
\end{align}
In the series
$ \sum_{r\ge 1} \frac{(t\Delta)^{r-1}}{r!} \Delta^{-1} (\partial^2_1- \partial^2_2)\partial_2 (u^0_1u^0_1) $,  the first term $ \Delta^{-1} (\partial^2_1- \partial^2_2)\partial_2 (u^0_1u^0_1)= -\partial_2 (u^0_1u^0_1) + 2 \partial^2_1\partial_2 (u^0_1u^0_1)$   contributes the main part and using the same way as Lemma \ref{NSlem2},  we have
\begin{align}
  t \left(\sum_{j\in \mathbb{N}_k}  2^{-jq} \left\|\triangle_j \partial_2 (u^0_1u^0_1)  \right\|^q_\infty \right)^{1/q}  \gtrsim  \eta  k^{1/q}.   \label{2d3.11}
\end{align}
Moreover,  taking notice of the definition of $a_l$ and $b_l$, one can repeat the procedures as in Lemma \ref{NSlem7}, Case 1 to obtain that
\begin{align}
  t \left(\sum_{j\in \mathbb{N}_k}  2^{-jq} \left\|\triangle_j \Delta^{-1}  \partial^2_1 \partial_2 (u^0_1u^0_1)  \right\|^q_\infty \right)^{1/q}  \lesssim \eta \varepsilon^2 k^{1/q}.   \label{2d3.12}
\end{align}
Similar to Lemma \ref{NSlem4}, we have
\begin{align}
  t \sum_{r\ge 2} \frac{1}{r!} \left( \sum_{j\in \mathbb{N}_k} 2^{-jq}\|   (t\Delta)^{r-1} \triangle_j     \partial_2 (u^0_1 u^0_1)  \|^q_\infty \right)^{1/q} \lesssim \eta     k^{1/q} 2^{-k}. \label{2d3.13}
\end{align}
Collecting the above estimates, we obtain that
$$\left(\sum_{j\in \mathbb{N}_k}  2^{-jq} \left\|\triangle_j \frac{\partial^2  u_1}{\partial\delta^2}\Big|_{\delta=0} \right\|^q_\infty \right)^{1/q} \gtrsim  \eta  k^{1/q}.$$
Recall that the estimates of higher order iterations in Section \ref{sect2.2} also hold for 2 dimensional case. So, the left part of the proof of Theorem \ref{NSIll} in 2D case proceeds in the same way as that of higher dimensional cases. The details of the proof are omitted.

\section{Other critical spaces } \label{other}

Even though NS is instable in the critical Besov spaces $\dot B^{-1}_{\infty,1}$, we still have some global well-posedness results in the subspaces of $\dot B^{-1}_{\infty,1}$.  Recall that Iwabuchi \cite{Iw10} (cf. also \cite{HuWa13, Wa06}) considered the well-posedness of NS in modulation spaces,   particularly in $M^{-1}_{\infty, 1}$ for which the norm is defined as follows.  Let $\sigma$ be a smooth
cut-off function with ${\rm supp} \ \sigma \subset  [-3/4, 3/4]^n$,
$\sigma_k=\sigma(\cdot-k)$ and
 $  \sum_{ k\in\mathbb{Z}^n}\sigma_k = 1.$
Denote
   $\square_k = \mathscr{F}^{-1} \sigma_k \mathscr{F},  \  k \in \mathbb{Z}^n. $ Put
\begin{align}
  \|f\|_{M^{-1}_{\infty,1}} = \sum_{k\in \mathbb{Z}^n} (1+|k|^2)^{-1/2} \|\Box_k f\|_\infty.  \label{3.1}
\end{align}
Modulation spaces was introduced by Feichtinger \cite{Fei83} (see also Gr\"ochenig \cite{Groch01}, here we adopt an equivalent norm; cf. \cite{WaHe07}). Lei and Lin \cite{LeLi11} considered the the global existence and uniqueness of 3D NS in the space $\mathcal{X}^{-1}$ for which the norm is defined by
\begin{align}
  \|f\|_{\mathcal{X}^{-1}} = \int_{\mathbb{R}^n} |\xi|^{-1} |\widehat{f}(\xi)| d\xi.   \label{3.2}
\end{align}

\begin{prop}\label{embed}
We have the following inclusions
\begin{align}
 \mathcal{X}^{-1} \subset  M^{-1}_{\infty,1} \subset B^{-1}_{\infty,1},   \label{3.3}
\end{align}
where $B^{-1}_{\infty,1}= L^p+ \dot B^{-1}_{\infty,1}$.
\end{prop}
{\bf Proof.} For the proof of the sharp embedding $ M^{-1}_{\infty,1} \subset B^{-1}_{\infty,1}$, one can referee to \cite{To04}, \cite{SuTo07}, \cite{WaHu07}. For any $f\in \mathcal{X}^{-1}$, by Young's inequality,
\begin{align}
\|f\|_{M^{-1}_{\infty,1}}  &  \lesssim \sum_{k\in \mathbb{Z}^n} \langle k\rangle^{-1} \|\sigma_k \widehat{f}\|_1   \nonumber\\
&    \lesssim  \sum_{k\in \mathbb{Z}^n}   \int_{k+[-1,1]^n} \langle \xi\rangle^{-1}   |\widehat{f}(\xi)| d\xi   \lesssim  \int_{\mathbb{R}^n}  |\xi|^{-1}   |\widehat{f}(\xi)| d\xi = \|f\|_{\mathcal{X}^{-1}}.   \label{3.4}
\end{align}
It follows that $\mathcal{X}^{-1} \subset  M^{-1}_{\infty,1}$. $\hfill \Box$

It is worth to mention that $M^{-1}_{\infty,1}$ and $\dot B^{-1}_{\infty,1}$ have no inclusions, since the lower (higher) frequency part of $\dot B^{-1}_{\infty,1}$ is smoother (rougher) than that of $M^{-1}_{\infty,1}$. However, we have $\mathcal{X}^{-1} \subset \dot B^{-1}_{\infty,1}$, which is easily seen by imitating the proof of \eqref{3.4}.
To some extent, $M^{-1}_{\infty,1}$ is also the critical space of NS, using similar ideas in this paper, we can show NS is ill-posed in $M^s_{\infty,1}$ if $s<-1$.

Let us observe the initial data $u^0_1$, we easily see that $\|u^0_1\|_{M^{-1}_{\infty,1}} \sim \|u^0_1\|_{\mathcal{X}^{-1}}\sim k$, $\|u^0_1\|_{\dot B^{-1}_{\infty,1}} \sim 1$. One sees that $ M^{-1}_{\infty,1} $ and $\mathcal{X}^{-1}$ are quite similar in the higher frequency part on which the norms of $u^0_1$ are far away from $B^{-1}_{\infty,1}$ or $\dot B^{-1}_{\infty,1}$ .

Finally, let us observe the initial data $u^0_1$ and the solution $u$ in critical Triebel-Lizorkin's spaces $\dot F^{-1}_{\infty, q}$ with $1<q<\infty$. Let us recall that   $\dot F^{-1}_{\infty, q}$ is defined by  (cf. \cite{Tr})
 \begin{align}
\dot F^{-1}_{\infty,q} = &  \Big\{ f\in \mathscr{Z}'(\mathbb{R}^n): \exists \ \{f_j\}_{j\in \mathbb{Z}} \ such \ that \ f = \sum_{j\in \mathbb{Z}} \triangle_j f_j   \nonumber\\
  &  \quad \quad \quad \quad  \quad \quad   in \ \mathscr{Z}'(\mathbb{R}^n) \ and \   \left\| \|\{2^{-j} f_j\}\|_{\ell^q}\right\|_\infty <\infty  \Big\}   \label{3.5}
\end{align}
for which the norm is
 \begin{align}
\|f\|_{\dot F^{-1}_{\infty,q}} = & \inf   \left\| \|\{2^{-j} f_j\}\|_{\ell^q}\right\|_\infty,  \label{3.6}
\end{align}
where the infimum is taken over all of the possible expressions in \eqref{3.5}.

Due to the embedding $\dot F^1_{1, q'} \subset \dot B^1_{1, q'}$, we immediately have from the duality that $\dot B^{-1}_{\infty, q} \subset \dot F^{-1}_{\infty, q}$.  Considering the initial data $u_0$ defined in \eqref{u0}, we have
 \begin{align}
\|u_0\|_{\dot F^{-1}_{\infty,q}}  \lesssim  \|u_0\|_{\dot B^{-1}_{\infty,q}} \lesssim 1, \ \ \forall q>1.  \label{3.7}
\end{align}
According to Koch and Tataru's result, NS has a unique small global solution in  $BMO^{-1}=\dot F^{-1}_{\infty,2}$. One may ask how does the solution vanish in $BMO^{-1}$ when $\delta \to 0$?  Now we outline the essential difference between $BMO^{-1}$ and  $\dot B^{-1}_{\infty,2}$ for NS with initial data as in \eqref{u0}.   By Lemma \ref{NSlem2}, we see that $(\partial_1-\partial_2)(u^0_1u^0_1)$ contributes the main part in critical Besov space  $\dot B^{-1}_{\infty,q}$ which grows like $O(k^{1/q})$. However, we can show it is very small in $\dot F^{-1}_{\infty,q}$  and we have
 \begin{align}
t\| (\partial_1-\partial_2)(u^0_1u^0_1)\|_{\dot F^{-1}_{\infty,q}} \lesssim  \varepsilon, \ \ t \le 2^{-2k}.  \label{3.8}
\end{align}
Taking $\tilde{\triangle}_j= \triangle_{j-1}+ \triangle_{j}+ \triangle_{j+1}$, we see that $f=\sum_j  \triangle_j (\tilde{\triangle}_j f)$ and
\begin{align}
t\| (\partial_1-\partial_2)(u^0_1u^0_1)\|_{\dot F^{-1}_{\infty,q}} & \le  \varepsilon t\left\| \sum_{j} |\tilde{\triangle}_j (u^0_1u^0_1)| \right\|_{\infty} \nonumber\\
 & \quad  +  t \left\| \sum_{j} |2^{-j}\tilde{\triangle}_j (\partial_1-\partial_2- \varepsilon 2^j)(u^0_1u^0_1)| \right\|_{\infty}. \label{3.9}
\end{align}
In view of \eqref{2.30}, we see that
$$\sum_j |\tilde{\triangle}_j (u^0_1u^0_1)| = 2^{2k+1} \sum_j |\Phi^+_j \Phi^-_j | +... \lesssim 2^{2k+1} \sum_j |\check{\varrho}(\cdot+a_j)|^2+... \lesssim 2^{2k},$$
which implies that the first term  in \eqref{3.9} is less than $\varepsilon$ as $t\le 2^{-2k}$, however, it is $O(k^{1/q})$ in critical Besov space $\dot B^{-1}_{\infty,q}$. The second term  in \eqref{3.9} is a remainder term which is exponentially small like $2^{- ck}$.  Hence, we have \eqref{3.8}.
Roughly speaking, for the solution $u$ considered in this paper,  the velocities localized in different frequencies, say $\triangle_j u$ and $\triangle_{m} u$ ($m\neq j$)  have concentrations in different physical regions which are far away from each other, so their superpositions in physical spaces have no inflation phenomena, this is why $u$ is small in $BMO^{-1}$.  \\

  \noindent{\bf Acknowledgment.} The author  is
supported in part by the National Science Foundation of China, grants   11271023.

\end{document}